\documentclass[11pt]{article}
\usepackage[utf8]{inputenc}
\usepackage{graphicx}
\usepackage{amsmath, amsthm, amssymb} 
\usepackage{hyperref} 
\usepackage{multicol}
\usepackage{multirow}
\usepackage{caption}
\usepackage{xcolor}
\usepackage{subcaption}
\usepackage[labelfont=bf]{caption}
\usepackage{pdflscape}
\usepackage[noend]{algpseudocode}

\usepackage{algorithm}

\newcommand{\myalgsheader}[0]

\algnewcommand{\IIf}[1]{\State\algorithmicif\ #1\ \algorithmicthen}
\algnewcommand{\EndIIf}{\unskip\ \algorithmicend\ \algorithmicif}
\algnewcommand{\IElse}[1]{\State\algorithmicelse\ #1\ \algorithmicthen}
\algnewcommand{\IfThenElse}[3]{% \IfThenElse{<if>}{<then>}{<else>}
  \State \algorithmicif\ #1\ \algorithmicthen\ #2\ \algorithmicelse\ #3}

\begin{document}
\title{Using Mixed Precision in Low-Synchronization Reorthogonalized Block Classical Gram-Schmidt\footnote{We acknowledge funding from the Charles University PRIMUS project
No. PRIMUS/19/SCI/11, Charles University GA UK project No. 202722, Charles University Research Program No. UNCE/SCI/023, and the Exascale Computing Project (17-SC-20-SC), a collaborative effort of the U.S. Department of Energy Office of Science and the National Nuclear Security Administration.}}

\author{Eda Oktay and Erin Carson}
\date{}
\maketitle

\paragraph{Abstract.}
%Mixed precision hardware has recently become commercially available, and more than 25\% of the supercomputers in the TOP500 list now have mixed precision capabilities. 
Using lower precision in algorithms can be beneficial in terms of reducing both computation and communication costs. 
%There are many current efforts towards developing mixed precision numerical linear algebra algorithms, which can lead to speedups in real applications. 
Motivated by this, we aim to further the state-of-the-art in developing and analyzing mixed precision variants of iterative methods. 
%
%In iterative methods based on Krylov subspaces, the orthogonal basis is generated by Arnoldi or Lanczos methods or their variants. In long recurrence methods such as GMRES, one needs to use an explicit orthogonalization scheme such as Gram-Schmidt to orthonormalize the vectors generated. 
%One of the variants of the Gram-Schmidt method is Classical Gram-Schmidt (CGS). From a performance perspective, CGS is preferable as it only involves one synchronization per vector. 
%
%Krylov subspace methods such as GMRES are typically communication-bound on modern machines; the runtime is usually dominated by the cost of global synchronizations which are required for the necessary orthogonalization. This has motivated the development of various algorithmic variants which attempt to reduce the communication cost while maintaining a stable algorithm. Recent work has focused on the development of low-synchronization block variants of Gram-Schmidt orthogonalization, which, when used within communication-avoiding ($s$-step) or block variants of GMRES, can reduce the number of global synchronizations to one per block. 
%
In this work, we focus on the block variant of low-synchronization classical Gram-Schmidt with reorthogonalization, which we call BCGSI+LS. We demonstrate that the loss of orthogonality produced by this orthogonalization scheme can exceed $O(u)\kappa(\mathcal{X})$, where $u$ is the unit roundoff and $\kappa(\mathcal{X})$ is the condition number of the matrix to be orthogonalized, and thus we can not in general expect this to result in a backward stable block GMRES implementation. We then develop a mixed precision variant of this algorithm, called BCGSI+LS-MP, which uses higher precision in certain parts of the computation. We demonstrate experimentally that for a number of challenging test problems, our mixed precision variant successfully maintains a loss of orthogonality below $O(u)\kappa(\mathcal{X})$. This indicates that we can achieve a backward stable block GMRES algorithm that requires only one synchronization per iteration.

\section{Introduction and Background} \label{sec:intro}
Block Gram-Schmidt (BGS) algorithms are used for computing the QR factorization of a given matrix $\mathcal{X}\in \mathbb{R}^{m\times n}$, $\mathcal{X}=\mathcal{Q}\mathcal{R}$, where $\mathcal{Q}$ is an orthogonal matrix and $\mathcal{R}$ is an upper triangular matrix. We assume $\mathcal{X}$ is partitioned into blocks of size $m\times s$, i.e., $\mathcal{X}=[X_1,\ldots,X_p]$, where $X_i\in \mathbb{R}^{m\times s}$, under the assumption of $n/s=p$. The matrices $\mathcal{Q}$ and $\mathcal{R}$ are partitioned in a similar manner, i.e., $\mathcal{Q}=[Q_1,\ldots,Q_p]$, where $Q_i \in \mathbb{R}^{m\times s}$ and $\mathcal{R}_{j,k} = R_{j,k}$, where $R_{j,k}\in\mathbb{R}^{s\times s}$. BGS algorithms are widely used in block Krylov subspace methods such as block GMRES \cite{bdj:06}, and in communication-avoiding Krylov subspace methods such as CA-GMRES \cite{b:14}. Using a block approach in Krylov subspace methods can improve performance by enabling the use of BLAS-3 operations, improving convergence behavior, and/or reducing the communication cost. 

In parallel settings, the inner product and norm operations within BGS algorithms require global reductions, i.e., synchronizations. This means that all compute nodes of a machine need to synchronize and exchange information. In large-scale settings, this can become a major computational bottleneck depending on the particular BGS variant and particular parallel setting. BGS algorithms also suffer from loss of orthogonality, i.e., the matrix $\mathcal{Q}$ is not exactly orthogonal when computed in finite precision. The loss of orthogonality is defined as the quantity $\|I-\bar{\mathcal{Q}}^T\bar{\mathcal{Q}}\|_2$, where $\bar{\mathcal{Q}}$ is the computed Q factor of $\mathcal{X}$. The loss of orthogonality plays an important role in the stability of Krylov subspace methods. For instance, it was proved in \cite{prs:06} that the level of orthogonality provided by the modified Gram-Schmidt (MGS) algorithm, i.e., $\|I-\bar{\mathcal{Q}}^T\bar{\mathcal{Q}}\|_2<\mathcal{O}(u)\kappa(\mathcal{X})$, where $u$ is the unit round-off (for double precision, $u\approx10^{-16}$), is sufficient to guarantee a backward stable GMRES algorithm. This motivates us to find a BGS variant that achieves at least $\mathcal{O}(u)\kappa(\mathcal{X})$ loss of orthogonality and simultaneously requires a small number of global synchronizations.  

One of the most common BGS variants is the block classical Gram-Schmidt (BCGS) algorithm \cite{s:03} given in Algorithm \ref{alg:bcgs}. Like all BGS variants, BCGS uses a non-block orthogonalization algorithm referred to as \verb*|IntraOrtho| for intrablock orthogonalization, such as Householder QR, CGS, MGS, or Cholesky QR. Note that BCGS requires one synchronization in line \ref{bcgs1} and potentially one or more synchronizations in line \ref{bcgs2}, depending on what is used as the \verb*|IntraOrtho|.

\begin{algorithm}[htbp!]
	\caption{BCGS \label{alg:bcgs}\cite{s:03}} 
	\begin{algorithmic}[1]
		\State{$[Q_1,R_{11}]$ = \verb*|IntraOrtho|($X_1$)}
		\For{$k = 1$: $p-1$}  
		\State{$R_{1:k,k+1}=Q_{1:k}^TX_{k+1}$} \label{bcgs1}
		\State{$W = X_{k+1}-Q_{1:k}R_{1:k,k+1}$} 
		\State{$[Q_{k+1},R_{k+1,k+1}]$= \verb*|IntraOrtho|($W$)} \label{bcgs2}
		\EndFor
	\end{algorithmic}
\end{algorithm}

According to the conjecture in \cite{k:74}, the loss of orthogonality in BCGS is bounded by $\|I-\bar{\mathcal{Q}}^T\bar{\mathcal{Q}}\|_2<\mathcal{O}(u)\kappa^{n-1}(\mathcal{X})$ as long as $\mathcal{O}(u)\kappa(\mathcal{X})<1$. Moreover, if the diagonal blocks of $\mathcal{R}$ are computed in a special way using Cholesky, the authors in \cite{clr:21} proved that $\|I-\bar{\mathcal{Q}}^T\bar{\mathcal{Q}}\|_2<\mathcal{O}(u)\kappa^2(\mathcal{X})$ holds provided that $\mathcal{O}(u)\kappa^2(\mathcal{X})<1$.

\subsection{Reorthogonalized Variants (BCGSI+)}
Round-off errors and cancellation causes BCGS to lose orthogonality which makes the algorithm unstable. To overcome this problem, one can use reorthogonalization.

To use reorthogonalization in the BCGS algorithm, instead of calculating the final $R$ in lines \ref{bcgs1}-\ref{bcgs2} of Algorithm \ref{alg:bcgs}, the orthogonalization is performed two times and the resulting $R$ factors are combined. This variant of BCGS (BCGSI+) given in Algorithm \ref{alg:bcgs+} was analyzed by the authors in \cite{bs:13}. The algorithm orthogonalizes for the first time in lines \ref{purple1}-\ref{purple2}. Then using the previously computed $\hat{Q}$, the vectors are orthogonalized for the second time in lines \ref{green1}-\ref{green2}. The $R$ factors are then combined in lines \ref{yellow1}-\ref{yellow2}. Note that BCGSI+ requires twice as many synchronizations as BCGS. 

\begin{algorithm}[htbp!]
	\caption{BCGSI+ \cite{bs:13} \label{alg:bcgs+}} 
	\begin{algorithmic}[1]
		\State{Allocate memory for $Q$ and $R$}
		\State{$[Q_1,R_{11}]$ = \verb*|IntraOrtho|($X_1$)} 
		\For{$k = 1$: $p-1$}  
		\State{$R_{1:k,k+1}^{(1)}=Q_{1:k}^TX_{k+1}$}\label{purple1}
		\State{$W = X_{k+1}-Q_{1:k}R_{1:k,k+1}^{(1)}$}
		\State{$[\hat{Q},R_{k+1,k+1}^{(1)}]$= \verb*|IntraOrtho|($W$)} \label{purple2}
		\State{$R_{1:k,k+1}^{(2)}=Q_{1:k}^T\hat{Q}$}\label{green1}
		\State{$W = \hat{Q}-Q_{1:k}R_{1:k,k+1}^{(2)}$}
		\State{$[Q_{k+1},R_{k+1,k+1}^{(2)}]$= \verb*|IntraOrtho|($W$)}\label{green2}
		\State{$R_{1:k,k+1}=R_{1:k,k+1}^{(1)}+R_{1:k,k+1}^{(2)}R_{k+1,k+1}^{(1)}$}\label{yellow1}
		\State{$R_{k+1,k+1}=R_{k+1,k+1}^{(2)}R_{k+1,k+1}^{(1)}$}\label{yellow2}
		\EndFor
	\end{algorithmic}
\end{algorithm}

The authors in \cite{bs:13} proved that if a method with $\|I-\bar{\mathcal{Q}}^T\bar{\mathcal{Q}}\|_2\leq\mathcal{O}(u)$ is used as \verb*|IntraOrtho| and $\mathcal{O}(u)\kappa(\mathcal{X})<1$, then BCGSI+ has a loss of orthogonality bounded by $\|I-\bar{\mathcal{Q}}^T\bar{\mathcal{Q}}\|_2\leq\mathcal{O}(u)$. 

\subsection{Low-Synchronization Variants (BCGSI+LS)}
%In parallel settings, the inner product and norm operations within BGS algorithms require global reductions, i.e., synchronizations. This means that all compute nodes of a machine need to synchronize and exchange information. In large-scale settings, this can become a major computational bottleneck. 
The goal of reducing the number of synchronizations required in Gram-Schmidt algorithms motivated work into developing so-called ``low-synchronization'' variants of BGS and other orthogonalization routines, which require only a single synchronization per block; see, e.g., \cite{y:20}. 

The low-sync BCGSI+ algorithm (BCGSI+LS) given in Algorithm \ref{alg:bcgsi+ls} was recently introduced in \cite{y:20} for use within GMRES. The BCGSI+LS algorithm is a block generalization of the CGSI+LS algorithm in \cite{s:21}, which is based on computing a strictly lower triangular matrix one block of rows at a time in a single global reduction, lagging the normalization step, and merging it into this single reduction. The development of this approach was based on the work of Ruhe \cite{r:83}; the authors of \cite{s:21} observed that MGS/CGS could be interpreted as a variant of Gauss-Seidel/Gauss-Jacobi iterations for solving the normal equations where the associated orthogonal projector is given as 
\[
I-Q_{1:k-1}T_{1:k-1,1:k-1}Q^T_{1:k-1}, \qquad \text{for}\quad T_{1:k-1,1:k-1}\approx (Q_{1:k-1}^T Q_{1:k-1})^{-1}.
\] 
For CGSI+LS, this $T$ is computed iteratively via $T_{1:k-1,1:k-1}\approx I-L_{1:k-1,1:k-1}-L^T_{1:k-1,1:k-1}$. According to \cite{c:22}, the algorithm can also be thought of as splitting $T_{1:k-1,1:k-1}$ into two parts, $I-L_{1:k-1,1:k-1}$ and a delayed reorthogonalization step $R_{1:k-1,k-1}=R_{1:k-2,k-1}-L^T_{k-1,1:k-2}$, the latter of which is applied in the next iteration. 

The block generalization of this idea leads to BCGSI+LS. We note in Algorithm \ref{alg:bcgsi+ls} that the block analogs of $T$ and $L$ above are not explicitly computed. We also note that the resulting BCGSI+LS algorithm has no explicit \verb*|IntraOrtho|, and lines \ref{1a}, \ref{2b}, and \ref{7a} are computed via Cholesky factorization. Asymptotically, BCGSI+LS has only one synchronization point per block, which occurs in lines \ref{1a} and \ref{2a} in Algorithm \ref{alg:bcgsi+ls}.

\begin{algorithm}[htbp!]
	\caption{BCGSI+LS \cite{y:20} \label{alg:bcgsi+ls}} 
	\newsavebox\aamatrix
	\savebox{\aamatrix}{$\begin{bmatrix} W & Z \\ \Omega& Y \end{bmatrix}$}%
	\newsavebox\ddmatrix
	\savebox{\ddmatrix}{$\begin{bmatrix} R^T_{k-1,k-1}R_{k-1,k-1}&P \end{bmatrix}$}%
	\newsavebox\ccmatrix
	\savebox{\ccmatrix}{$\begin{bmatrix} U &X_k \end{bmatrix}$}%
	\newsavebox\eematrix
	\savebox{\eematrix}{$\begin{bmatrix} \mathcal{Q}_{1:k-2} & U \end{bmatrix}$}%
	\newsavebox\ffmatrix
	\savebox{\ffmatrix}{$\begin{bmatrix} U & X_k \end{bmatrix}$}%\\
	\newsavebox\ggmatrix
	\savebox{\ggmatrix}{$\begin{bmatrix} W & Z \\ \Omega& Z \end{bmatrix}$}%
	\newsavebox\iimatrix
	\savebox{\iimatrix}{$\begin{bmatrix} \Omega & Y \end{bmatrix}$}%
	\newsavebox\jjmatrix
	\savebox{\jjmatrix}{$\begin{bmatrix} W &Z \end{bmatrix}$}%
	\newsavebox\kkmatrix
	\savebox{\kkmatrix}{$\begin{bmatrix} W  \\ \Omega \end{bmatrix}$}%
	\newsavebox\llmatrix
	\savebox{\llmatrix}{$\begin{bmatrix} \mathcal{Q}_{1:s-1} & U \end{bmatrix}$}%
	\begin{algorithmic}[1]
		\State{Allocate memory for $\mathcal{Q}$ and $\mathcal{R}$}
		\State{$U = X_1$} 
		\For{$k=2,\ldots,p$}  
		\If{$k=2$} 
		\State $\usebox{\ddmatrix}=U^T\usebox{\ccmatrix}$ \label{1a}
		\ElsIf{$k>2$} 
		\State {$ \usebox{\aamatrix} = \usebox{\eematrix}^T\usebox{\ffmatrix}$}\label{2a}
		\State $\usebox{\ddmatrix}=\usebox{\iimatrix}-W^T\usebox{\jjmatrix}$ \label{2b}
		\EndIf
		\State{$R_{k-1,k}=R^{-T}_{k-1,k-1}P$}\label{3a}
		\If{$k=2$} 
		\State{$Q_{k-1}=UR^{-1}_{k-1,k-1}$} \label{4a}
		\ElsIf{$k>2$} 
		\State{$\mathcal{R}_{1:k-2,k-1}=\mathcal{R}_{1:k-2,k-1}+W$}
		\State{$\mathcal{R}_{1:k-2,k}=Z$}
		\State{$Q_{k-1}=(U-\mathcal{Q}_{1:k-2}W)R^{-1}_{k-1,k-1}$}\label{5a}
		\EndIf
		\State{$U=X_k-\mathcal{Q}_{1:k-1}\mathcal{R}_{1:k-1,k}$}
		\EndFor
		\State{$\usebox{\kkmatrix} = \usebox{\llmatrix}^TU$}\label{6a}
		\State{$R^T_{s,s}R_{s,s}=\Omega - W^TW$}\label{7a}
		\State{$\mathcal{R}_{1:s-1,s}=\mathcal{R}_{1:s-1,s}+W$}
		\State{$Q_s=(U-\mathcal{Q}_{1:s-1}W)R^{-1}_{s,s}$}\label{8a}
		\State \Return{$\mathcal{Q}=[Q_1,\ldots,Q_s], \mathcal{R}=(R_{jk})$}
	\end{algorithmic}
\end{algorithm}

A conjecture in \cite{c:22} states that if $\mathcal{O}(u)\kappa^2(\mathcal{X})<1$, then the loss of orthogonality of BCGSI+LS satisfies $\|I-\bar{Q}^T\bar{Q}\|_2<\mathcal{O}(u)\kappa^2(\mathcal{X})$. Thus although BCGSI+LS has the performance advantage that it only requires one synchronization per block, such a significant loss of orthogonality may make it unsuitable for use within GMRES. This motivates us to try to selectively use higher precision in some parts of the BCGSI+LS algorithm to decrease the loss of orthogonality while still maintaining a single synchronization per block. 

\section{Mixed Precision BCGSI+LS (BCGSI+LS-MP)}
%Using low precision in an algorithm can save communication and computation costs. For instance, it is observed that using half-precision (float 16) is $16\times$ faster than using single precision (float 32) in NVIDIA A100 GPU \cite{a100}. However, using low precision throughout an algorithm may not be sufficient for the algorithms requiring high accuracy. This motivates the search for mixed precision algorithms, which selectively use high and low precision in different parts of the computation to achieve both performance and accuracy.
Our mixed precision approach, which we call BCGSI+LS-MP, is given in Algorithm \ref{alg:bcgsi+ls_mp}. For a working precision $u$, we use the higher precision $u^2$ in two aspects of the algorithm: for computing the Cholesky factorizations in lines \ref{1aa}, \ref{2aa}/\ref{2ab}, and \ref{6aa}/\ref{7aa}, and in applying the corresponding inverses of the $R$ factors in lines \ref{3aa}, \ref{4aa}, \ref{5aa}, and \ref{8aa}. 

We note that BCGSI+LS-MP still has only one synchronization point per block, which occurs in lines \ref{1aa} and \ref{2aa}. Since we now use precision $u^2$ in these lines, this means that we have doubled the size of the reduction; i.e., we have doubled the bandwidth and computation costs. In the latency-bound regime, where low-synchronization algorithms are most beneficial, this overhead may not be significant, in particular, in cases where precisions $u$ and $u^2$ are both implemented in hardware. The higher precision computations in lines \ref{2ab}, \ref{3aa}, \ref{4aa}, and \ref{5aa} are all local computations, and thus we expect the extra overhead to be insignificant. We note that lines \ref{6aa}, \ref{7aa}, and \ref{8aa} are only computed once at the very end of the algorithm. Overall, the resulting overhead of using mixed precision will be highly dependent on the particular problem size and machine parameters; a performance study will be the subject of future work. 

%From our numerical experiments as shown in Section \ref{sec:exp}, when $u^2$ is used instead of $u$ in some lines in Algorithm \ref{alg:bcgsi+ls}, the upper bound for the loss of orthogonality of the mixed precision BCGSI+LS algorithm (BCGSI+LS-MP) becomes $\mathcal{O}(u)\kappa(\mathcal{X})$, resulting in the possibility of backward stable solution when used with the GMRES algorithm. Furthermore, we see that using the mixed precision idea, BCGSI+LS-MP may be applied to worse-conditioned matrices than BCGSI+LS. The BCGSI+LS-MP algorithm is given in Algorithm \ref{alg:bcgsi+ls_mp}.

\begin{algorithm}[htbp!]
	\caption{BCGSI+LS-MP \label{alg:bcgsi+ls_mp}} 
	\newsavebox\aaamatrix
	\savebox{\aaamatrix}{$\begin{bmatrix} W & Z \\ \Omega& Y \end{bmatrix}$}%
	\newsavebox\dddmatrix
	\savebox{\dddmatrix}{$\begin{bmatrix} R^T_{k-1,k-1}R_{k-1,k-1}&P \end{bmatrix}$}%
	\newsavebox\cccmatrix
	\savebox{\cccmatrix}{$\begin{bmatrix} U &X_k \end{bmatrix}$}%
	\newsavebox\eeematrix
	\savebox{\eeematrix}{$\begin{bmatrix} \mathcal{Q}_{1:k-2} & U \end{bmatrix}$}%
	\newsavebox\fffmatrix
	\savebox{\fffmatrix}{$\begin{bmatrix} U & X_k \end{bmatrix}$}%\\
	\newsavebox\gggmatrix
	\savebox{\gggmatrix}{$\begin{bmatrix} W & Z \\ \Omega& Z \end{bmatrix}$}%
	\newsavebox\iiimatrix
	\savebox{\iiimatrix}{$\begin{bmatrix} \Omega & Y \end{bmatrix}$}%
	\newsavebox\jjjmatrix
	\savebox{\jjjmatrix}{$\begin{bmatrix} W &Z \end{bmatrix}$}%
	\newsavebox\kkkmatrix
	\savebox{\kkkmatrix}{$\begin{bmatrix} W  \\ \Omega \end{bmatrix}$}%
	\newsavebox\lllmatrix
	\savebox{\lllmatrix}{$\begin{bmatrix} \mathcal{Q}_{1:s-1} & U \end{bmatrix}$}%
	\begin{algorithmic}[1]
		\State{Allocate memory for $\mathcal{Q}$ and $\mathcal{R}$}
		\State{$U = X_1$} 
		\For{$k=2,\ldots,p$}  
		\If{$k=2$} 
		\State $\usebox{\dddmatrix}=U^T\usebox{\cccmatrix}$ {\hspace*{\fill} in precision $u^2$}\label{1aa}
		\ElsIf{$k>2$} 
		\State {$ \usebox{\aaamatrix} = \usebox{\eeematrix}^T\usebox{\fffmatrix}$} {\hspace*{\fill} in precision $u^2$}\label{2aa}
		\State $\usebox{\dddmatrix}=\usebox{\iiimatrix}-W^T\usebox{\jjjmatrix}$ {\hspace*{\fill} in precision $u^2$}\label{2ab}
		\EndIf
		\State{$R_{k-1,k}=R^{-T}_{k-1,k-1}P$}{\hspace*{\fill} in precision $u^2$}\label{3aa}
		\If{$k=2$} 
		\State{$Q_{k-1}=UR^{-1}_{k-1,k-1}$}{\hspace*{\fill} in precision $u^2$} \label{4aa}
		\ElsIf{$k>2$} 
		\State{$\mathcal{R}_{1:k-2,k-1}=\mathcal{R}_{1:k-2,k-1}+W$}{\hspace*{\fill} in precision $u\phantom{^2}$}
		\State{$\mathcal{R}_{1:k-2,k}=Z$}{\hspace*{\fill} in precision $u\phantom{^2}$}
		\State{$Q_{k-1}=(U-\mathcal{Q}_{1:k-2}W)R^{-1}_{k-1,k-1}$} {\hspace*{\fill} in precision $u^2$}\label{5aa}
		\EndIf
		\State{$U=X_k-\mathcal{Q}_{1:k-1}\mathcal{R}_{1:k-1,k}$}{\hspace*{\fill} in precision $u\phantom{^2}$}
		\EndFor
		\State{$\usebox{\kkkmatrix} = \usebox{\lllmatrix}^TU$} {\hspace*{\fill} in precision $u^2$}\label{6aa}
		\State{$R^T_{s,s}R_{s,s}=\Omega - W^TW$} {\hspace*{\fill} in precision $u^2$}\label{7aa}
		\State{$\mathcal{R}_{1:s-1,s}=\mathcal{R}_{1:s-1,s}+W$}{\hspace*{\fill} in precision $u\phantom{^2}$}
		\State{$Q_s=(U-\mathcal{Q}_{1:s-1}W)R^{-1}_{s,s}$} {\hspace*{\fill} in precision $u^2$}\label{8aa}
		\State \Return{$\mathcal{Q}=[Q_1,\ldots,Q_s], \mathcal{R}=(R_{jk})$}
	\end{algorithmic}
\end{algorithm}

%However, it is essential to note that there is a trade-off between the loss of orthogonality and computational cost. Due to the use of float 128 in a large part of the algorithm, BCGSI+LS-MP may be more computationally expensive than BCGSI+LS. Although BCGSI+LS-MP increases the computation and bandwidth costs, since we are computing with and moving more bits, it keeps the number of messages the same. Therefore, the resulting overhead will depend on the particular latency and bandwidth parameters of the target machine.

\section{Numerical Experiments}\label{sec:exp}
We now seek to demonstrate numerically, on a set of challenging test problems, that our mixed precision approach BCGSI+LS-MP improves the loss of orthogonality relative to the uniform precision approach BCGSI+LS. 

To illustrate the comparison of the methods in terms of the loss of orthogonality, we performed numerical experiments in MATLAB using the block Gram-Schmidt variants available at \verb*|github.com/katlund/BlockStab| with L\"{a}uchli, monomial, and glued matrices. Each of the matrices has dimensions $[m, p, s]$, where $m$ is the number of rows, $p$ is the number of block vectors, and $s$ is the number of columns per block vector. The widths of blocks are specified by the input \verb*|svec|. For a detailed investigation of these matrices on BGS variants, see \cite{c:22}. The experiments are performed on a computer with AMD Ryzen 5 4500U having 6 CPUs and 8 GB RAM with OS system Ubuntu 22.04 LTS. In our numerical experiments, we used double precision in MATLAB for the working precision $u$, and quadruple precision for $u^2$. We used the Advanpix toolbox \cite{advanpix} to simulate quadruple precision.

Each plot shows the loss of orthogonality versus condition number for matrices $\mathcal{X}$ of a given type. The dashed black line represents the $\mathcal{O}(u)\kappa(\mathcal{X})$ bound and the black solid line represents the $\mathcal{O}(u)\kappa^2(\mathcal{X})$ bound. The dashed black line thus represents the loss of orthogonality bound for MGS, which is notable since MGS-GMRES is known to be backward stable \cite{prs:06}. We can thus conjecture that under certain constraints on the input matrix, an orthogonality scheme that provides this level of orthogonality can be expected to result in a backward stable GMRES implementation. The algorithm following the $\circ$ notation in the legends indicates the algorithm that is used as the \verb*|IntraOrtho| (which in our experiments, is always Householder QR). 
%{\color{red}If a method provides results below the dashed line, then the resulting GMRES algorithm can be guaranteed to be backward stable. However, if it is above this bound, GMRES may not provide a backward stable solution.}

\begin{figure}[!ht]
	\begin{minipage}{80mm}
		\includegraphics[width=\linewidth,trim={3cm 8cm 3cm 8cm},clip]{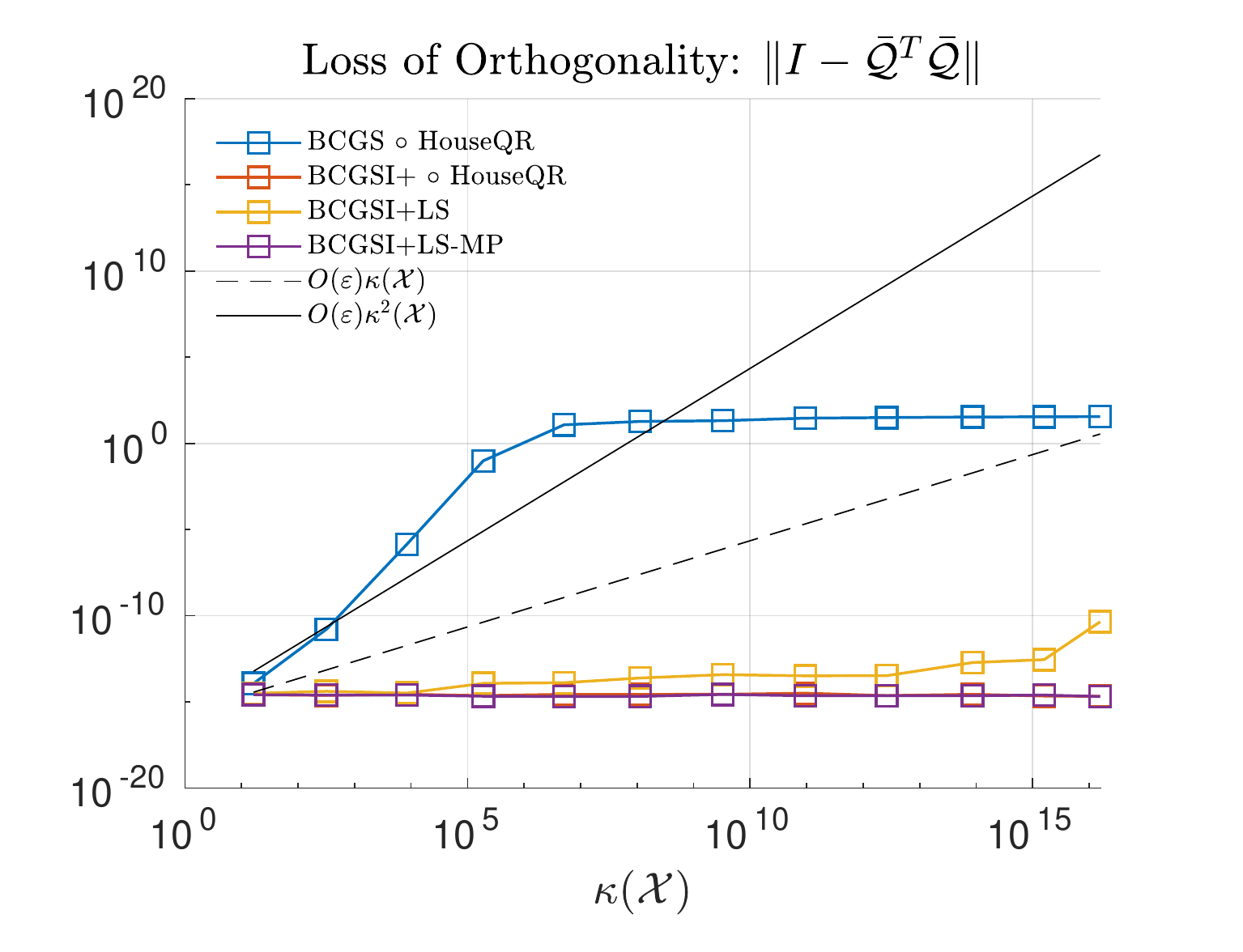}
		\caption{The loss of orthogonality of BCGS (blue), BCGSI+ (red), BCGSI+LS (yellow), and BCGSI+LS-MP (purple) versus condition number for glued matrices.}
		\label{fig:glued}
	\end{minipage}
	\hfill
	\begin{minipage}{80mm}
		\vspace*{-0.2cm}
		\includegraphics[width=\linewidth,trim={3cm 8cm 3cm 8cm},clip]{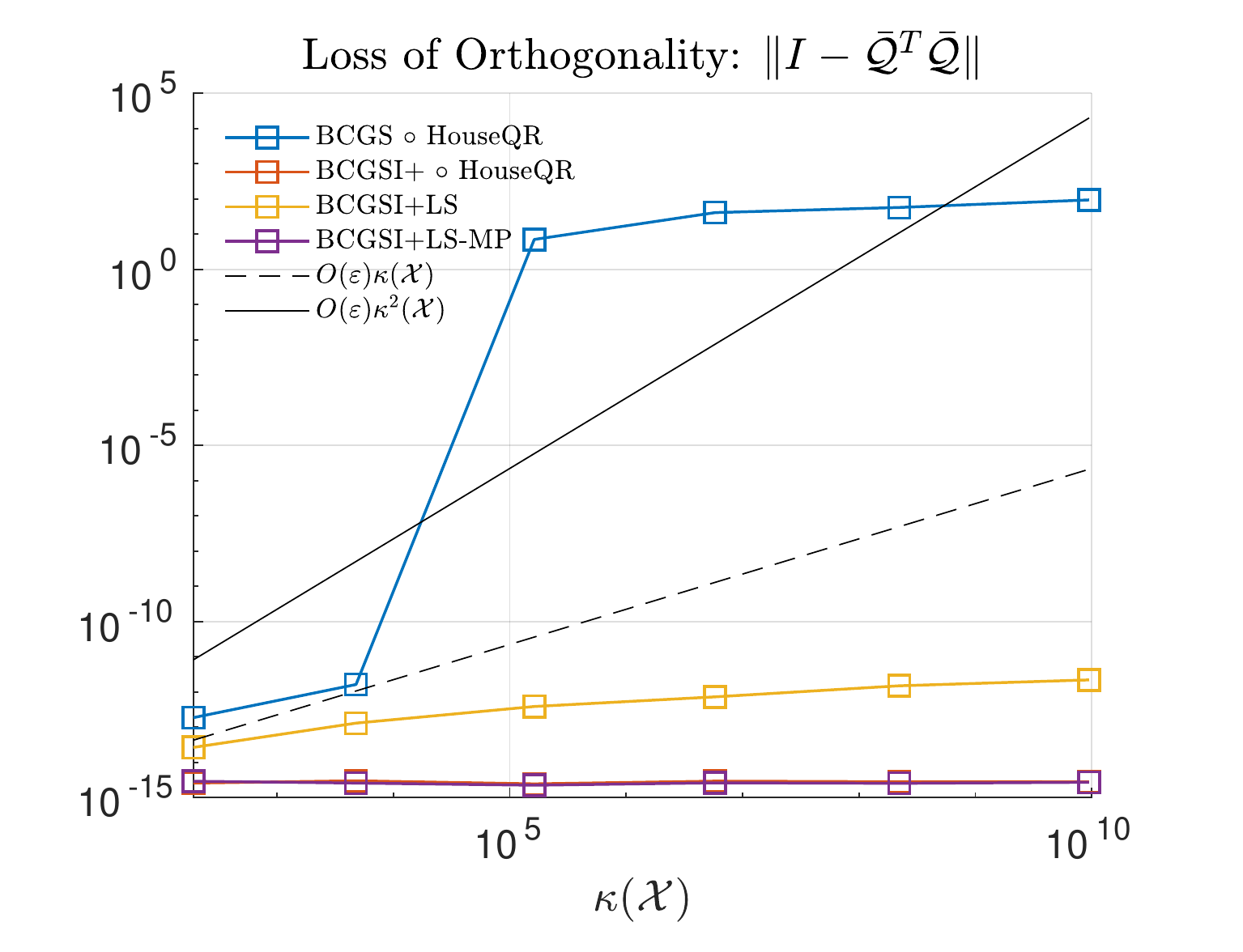}
		\caption{The loss of orthogonality of BCGS (blue), BCGSI+ (red), BCGSI+LS (yellow), and BCGSI+LS-MP (purple) versus condition number for monomial matrices.}
		\label{fig:monomial}
	\end{minipage}
	%	\hfill
	%	\begin{minipage}{80mm}
		%		\vspace{0.2cm}
		%		\includegraphics[width=\linewidth,trim={3cm 8cm 3cm 8cm},clip]{laeuchli}
		%		\caption{The loss of orthogonality of BCGS (blue), BCGSI+LS-MP (purple), BCGSI+LS (yellow), and BCGSI+ (red) when used with Householder QR as intrablock orthogonalization method for L\"{a}uchli matrix.}
		%		\label{fig:laeuchli}
		%	\end{minipage}
\end{figure}

The glued matrices introduced in \cite{sbl:06} are $m\times n$ matrices, where $n=nglued \times nbglued$, $nglued$ is the number of columns in a block, and $nbglued$ is the number of blocks that are glued together. For this study, we use glued matrices with dimension $[m, p, s]=[1000, 50, 4]$ with \verb*|svec=1:12|. From Figure \ref{fig:glued}, we see that even when an unconditionally stable intrablock orthogonalization method is used, the loss of orthogonality in BCGS can exceed the $\mathcal{O}(u)\kappa^2(\mathcal{X})$ bound. On the other hand, when reorthogonalization is used, the loss of orthogonality remains on the level $\mathcal{O}(u)$ (note that the red and purple markers overlap). Whereas the loss of orthogonality for BCGSI+LS starts to deviate from the level $\mathcal{O}(u)$ for larger condition numbers, for BCGSI+LS-MP it remains on the level $\mathcal{O}(u)$. 

The monomial test matrices are matrices $\mathcal{X}$ consisting of $p$ block vectors $X_k=[v_k|Av_k|\cdots|A^{s-1}v_k]$, $k=1,\ldots,p$, where $A$ is a diagonal $m\times m$ operator with evenly distributed eigenvalues in $(0.1,10)$, and $v_k$ are normalized randomly generated vectors from the uniform distribution. For this study, we use monomial matrices with dimension $[m, p, s]=[1000, 120, 2]$ with \verb*|svec=2:2:12|. We see from Figure \ref{fig:monomial} that the behavior of BCGS and reorthogonalized variants are the same as for the glued matrices: BCGS exceeds the $\mathcal{O}(u)\kappa^2(\mathcal{X})$ bound, and using reorthogonalization helps to decrease the loss of orthogonality to below $\mathcal{O}(u)\kappa(\mathcal{X})$. Also similarly to the glued matrices, the loss of orthogonality in BCGSI+LS deviates from $\mathcal{O}(u)$ (more significantly in this case), whereas for BCGSI+LS-MP, it remains on the level $\mathcal{O}(u)$. 

\begin{figure}[!ht]
	%\sidecaption
	\centering
	\includegraphics[width=85mm,trim={3cm 8cm 3cm 8cm},clip]{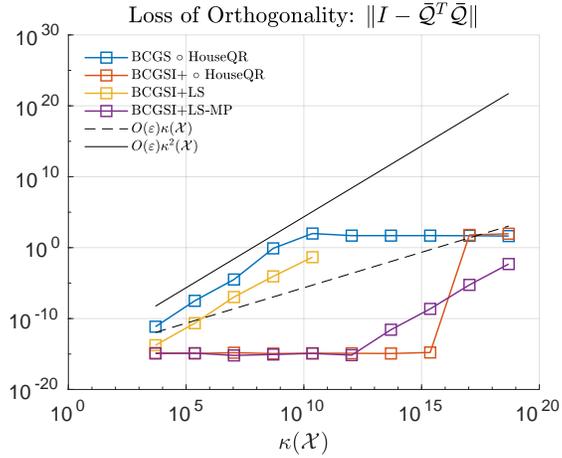}%
	\caption{The loss of orthogonality of BCGS (blue), BCGSI+ (red), BCGSI+LS (yellow), and BCGSI+LS-MP (purple) versus condition number for L\"{a}uchli matrices.}
	\label{fig:laeuchli}
\end{figure}

Our final, most interesting test case are L\"{a}uchli matrices \cite{l:61} of the form 
\[
\mathcal{X}= \begin{bmatrix}
	1 & 1& \ldots &1 \\
	\eta& & & \\
	&\eta& &\\
	& & \ddots & \\
	& & &\eta
\end{bmatrix}, \eta \in (\epsilon,\sqrt{\epsilon}),
\]
where $\eta$ is drawn randomly from a scaled uniform distribution. For L\"{a}uchli matrices, the columns are only barely independent, and the entries are either zero or close to $\epsilon$ which causes a high cancellation rate. Because of its structure, L\"{a}uchli matrices are often used in numerical experiments for illuminating the effects of finite precision error in Gram-Schmidt algorithms. For this study, we use L\"{a}uchli matrices with dimension $[m, p, s]=[1000, 100, 5]$ with \verb*|svec=logspace(-1,-16,10)|. Figure \ref{fig:laeuchli} shows the loss of orthogonality for each BCGS variant mentioned above. From the figure, we see that the loss of orthogonality of BCGS is above the $\mathcal{O}(u)\kappa(\mathcal{X})$ bound. When reorthogonalization is used, we observe that BCGSI+ stays on the level $\mathcal{O}(u)$ as expected. However, unlike in previous test cases, the loss of orthogonality in BCGSI+LS is now significantly worse. We see here that it exceeds the level $\mathcal{O}(u)\kappa(\mathcal{X})$, and the scaling behavior is on par with $\mathcal{O}(u)\kappa^2(\mathcal{X})$ (which is the bound that was conjectured in \cite{c:22}). The use of mixed precision in BCGSI+LS-MP seems to remedy this problem. The loss of orthogonality for BCGSI+LS-MP stays on the level $\mathcal{O}(u)$ at least until $\kappa(\mathcal{X}) \approx 10^{12}$.

\section{Conclusion and Discussion}
Orthogonalization processes are the core of Krylov subspace algorithms. For block Krylov subspace methods, block orthogonalization processes must be used to improve the performance. There are several block Gram-Schmidt variants, and each of these algorithms has different properties in terms of loss of orthogonality, communication, and computation costs. Recent work has focused on developing low-synchronized variants of (block) Gram-Schmidt algorithms, which require only a single global synchronization per (block) column \cite{y:20}. However, as our numerical experiments demonstrate, this reduced synchronization can come at the cost of decreased stability in terms of loss of orthogonality. 

We present a new block Gram-Schmidt variant called BCGSI+LS-MP, a variant of BCGSI+LS which uses higher precision in certain parts of the computation. Our numerical experiments demonstrate that this use of mixed precision can lead to a more stable algorithm than uniform precision BCGSI+LS while still requiring only a single global synchronization per block. We expect that in the latency-bound regime, where such low-synchronization algorithms are used, the increased bandwidth and computation costs due to the use of higher precision may be negligible. 

We note that here we have only provided empirical results; a rigorous proof of the loss of orthogonality in BCGSI+LS-MP remains future work. A first step towards this will be to prove a bound on the loss of orthogonality for uniform precision BCGSI+LS, which is currently missing in the literature. 
Furthermore, this preliminary study is based on the numerical experiments performed in MATLAB, which cannot give a good indication of resulting performance in large-scale parallel settings. A thorough performance study measuring the overhead of the use of higher precision within BCGSI+LS-MP is needed. 

\clearpage
\bibliographystyle{siamplain}
\bibliography{references}

\end{document}